# A novel, Fourier series based method of control optimization and its application to a discontinuous capsule drive model

*Sandra Zarychta, Tomasz Sagan, Marek Balcerzak[1], Artur Dabrowski, Andrzej Stefanski, Tomasz Kapitaniak*

*Division of Dynamics, Lodz University of Technology*

*Stefanowskiego 1/15, Lodz, Poland*

## Abstract

This paper presents a novel, Fourier series based numerical method of open-loop control optimization. Due to its flexible assumptions, it can be applied in a large variety of systems, including discontinuous ones or even "black boxes", whose equations are not fully known. This aspect is particularly important in mechanical systems, where friction or impact induced discontinuities are common. The paper includes a mathematical background of the new method, a detailed discussion of the algorithm and a numerical example, in which control function of a discontinuous capsule drive is optimized. It is expected that the proposed method can facilitate research in all areas where control of non-smooth, discontinuous or "black box" systems is crucial. In particular, authors hope that the presented algorithm is going to be used for control optimization of other capsule drives.

## 1. Introduction

A control of a dynamical system is optimal if it allows to perform a desired action at the lowest cost. Obviously, determination of such control is desirable in various fields of science and engineering [1, 2]. It is not only expected that devices, systems and processes surrounding us act as it was planned and provide expected results. It is also required that goals are achieved without wasting time, energy or other precious resources, i.e. optimally, according to some predefined criterion. Therefore, development of optimal control is an important issue not only from the purely scientific point of view. The question of optimal control is vital in any area, in which limited resources must be used efficiently, for example in engineering or in economics [1].

As this work focuses on optimal control of mechanical systems, assume that the control object is described by a set of ordinary differential equations (ODEs) [1]. In such case, the necessary conditions for optimality of a solution of an optimal control problem (OCP) have been described in the classical work by Pontryagin [3] in terms of variational calculus. These conditions are commonly referred to as Pontryagin's minimum principle (PMP) [2, 4, 5]. They are founded on the Lagrange multipliers approach [6] and described by means of so called Hamiltonian control function [1, 2], which attains its minimum over the whole trajectory of the system if its control is optimal.

Other classical result in the optimal control theory results from the dynamic programming method. Fundamental results in this field have been obtained by Bellman et al. [7, 8]. Dynamic programming is a general approach to optimization problems in which the task is divided into simpler subproblems in a recursive manner. In the context of optimal control, the dynamic programming approach is founded on a so-called cost-to-go function [2]. The cost-to-go function answers the following question: what is the minimal cost of transition of the system from the current state to the desired final state ? By means of such function it is possible to apply the dynamic programming philosophy of breaking the original problem into simpler subproblems. For instance, in order to minimize the cost over the whole trajectory, one can minimize the cost-to-go function from an intermediate state to the final state and from the initial state to the intermediate one, and thus divide control optimization over the time interval into cost minimization in two sub-intervals [1]. At the limit, as the length of the time interval approaches zero, the widely known the Hamilton-Jacobi-Bellman (HJB) partial differential equation is obtained. It is the sufficient condition for control optimality [1].

---

[1] Corresponding author, e-mail: marek.balcerzak.1@p.lodz.pl

In some simple problems [1, 2, 4], the Pontryagin's minimum principle (PMP) allows direct calculation of the function which meets the necessary conditions for optimality. However, in most cases, determination of optimal control requires application of numerical methods. The vast family of numerical algorithms for solving OCPs is divided into three main classes [4, 5].

- Dynamic programming,

- indirect methods based on calculus of variations,

- direct methods.

The first class encompasses methods based on results by Bellman et. al. [7, 8], particularly the HJB equation. In most cases, the cost-to-go function must be approximated numerically. Possible approaches include function expansion, finite difference or finite element techniques [4]. This method can be practically used only for optimal control problems (OCPs) in which dimension of the state space is low (except for the special case of the linear-quadratic regulator – LQR) [1]. The second class includes algorithms founded on the Pontryagin's Minimum Principle (PMP). The OCP is then reduced to solution of a boundary value problem [1], which usually requires application of numerical procedures. Exemplary methods include: single shooting (parametrization of the whole control function at once), multiple shooting (parametrization the control function in subsequent time intervals) and collocation methods (in which both state and control are parametrized using a predefined space of candidate solutions) [4, 5]. Note that the PMP yields necessary, but not sufficient, conditions of control optimality. The last class covers methods based on direct discretization of the OCP and thus transforming it into a nonlinear programming problem (NLP) to be optimized using a selected algorithm. Similarly as in the previous class, applicable numerical methods include single or multiple shooting methods, as well as collocation algorithms [4, 5].

Recently, an interesting tool named CASADI has been developed [9]. This library includes methods for nonlinear optimization and algorithmic differentiation, which allows effective estimation of derivatives of functions defined in terms of a computer program code. Therefore, it provides a complete framework allowing to solve multiple OCPs.

The OCP becomes more complicated when the system to be controlled is discontinuous, i.e. when the vector field, which defines the set of ODEs, is a discontinuous function. This often occurs in mechanical systems, in which impacts or dry friction play an important role [10, 11]. In such cases, the standard methods for solving OCPs fail. The PMP requires that the vector field is at least once continuously differentiable with respect to all its arguments [3], which precludes its application for non-smooth systems. Consequently, typical algorithms based on indirect numerical methods cannot be applied in such problems. Although the HJB equation does not involve differentiation of the vector field, so its smoothness is not formally required, it is reported that "optimize the discretization" strategy, being the foundation of numerical solutions of the HJB equations, causes excessive errors when vector field is discontinuous [11].

Several attempts have been made to solve OCP of non-smooth and discontinuous systems. In the works [12, 13] OCP of non-smooth, but continuous, systems has been discussed, whereas the publication [14] covers OCP of objects whose vector field is set-valued, but Lipschitz continuous. The paper [15] describes an approach based on computation of Jacobian matrix of the trajectory with respect to the initial conditions. This is possible as long as neither the initial state nor the final state is on the discontinuity. According to [15], this method leads to analog of PMP conditions. However, computation of this Jacobian matrix is not a simple task. Another possibility is application of neural networks and evolutionary computation [16]. Last but not least, the publication [11] describes results of approximating the non-smooth vector field with smooth functions, which offers promising results, but in fact a qualitatively different system is analyzed.

An interesting class of discontinuous control systems are so-called capsule drives. These devices, usually capsule-shaped, use an internal oscillator to produce inertia forces which, in the presence of dry friction, allow to move the capsule in a desired direction. Therefore, external moving parts, such as wheels, tracks, robotic arms etc., are no longer needed to produce motion. An interesting example is the vibro-impact drive [17], i.e. the drive in which a mass-on-spring oscillator is accompanied by a second spring, which can be impacted by the mass during its vibrations. Such arrangement causes that the resultant, horizontal force acting on the capsule is not symmetric and the capsule is able to move forward. Significant research on dynamics of such system and its modifications has already been conducted. The paper [17] provides a detailed bifurcation study. Optimization of progression of the capsule is presented, but only simple harmonic control functions are taken into account. The numerical studies have been verified experimentally in [18, 19]. The publications [20-22] present extensive analysis of different friction models and in [21] control functions which enable bidirectional motion are considered. The papers [23-25] contain a more detailed study of the harmonic control, including maximization of the rate of progression, as well as optimization of the energy consumption. Multistability control of the vibro-impact drive is available in [26]. Design and testing of a small-scale prototype are described in [27].

Other example is a pendulum capsule drive, in which the mass-on-spring oscillator is replaced by a pendulum. Its dynamics, analyzed in the paper [28], is somewhat more complex than in the previous case, as swinging of the pendulum influence contact force between the capsule and the underlying surface, which in turn affects the friction force. Different aspects of the system under consideration have already been investigated. The influence of viscoelasticity has been considered in [29] and the friction-induced hysteresis in [30]. Design and parameters optimization of a pre-designed control function profile is considered in the papers [31, 32]. In [33] the controller has been improved by implementation of a neural network based adaptation mechanism. Energy-related issues connected with the control function are studied in [34].

It seems that, in general, the research on control optimization of the capsule drives, which have been conducted up to this moment, focus on the following strategy. Firstly, select a "shape" or "profile" of a control function and parametrize it. The selected profile was usually sinusoidal in papers connected with the vibro-impact drive [17, 23-25], whereas an object-specific shape was commonly used in works concerning the pendulum capsule drive [31-33]. Then, optimize parameters of the profile. However, it is expected that interesting results could be obtained by a method, in which the shape of the control function is not assumed a priori, but is subject to optimization as well.

Such algorithm is presented in this paper. Authors would like to present a novel, Fourier series based numerical method of bounded optimal control estimation, which is intended to efficiently solve the optimal control problem (OCP) in non-smooth mechanical systems, including these with discontinuities resulting from dry friction or impacts. It is expected that the developed method is going to be applied in discontinuous systems of different kinds, as well as in objects whose dynamics is not fully known a priori ("black box" systems). In particular, it is shown that the algorithm is able to solve the problem of open-loop control optimization of capsule drives. In order to demonstrate correctness, efficiency and usefulness of the novel method, numerical simulation has been conducted, in which control of the pendulum capsule drive has been optimized.

## 2. Description of the method

In this chapter, the problem to be solved is precisely stated and the mathematical tools to be applied throughout its solution are introduced. Then, the proposed method of control optimization is precisely described and its correctness is proved.

### 2.1. Problem statement

Consider a controlled dynamical system described by the following ODE [1, 2]:

$$\dot{x}(t) = f[x(t), u(t), t], \qquad x(t_0) = x_0 \tag{1}$$

where $t \in \mathbb{R}$ is the time, $x(t) \in \mathbb{R}^n$ is a state vector at the time $t$, $\dot{x} = \frac{dx}{dt}$, $n$ is the order of the system, $u(t) \in \Omega \subset \mathbb{R}^r$ is a control vector (at the time $t$) of a dimension $r$, $f: \mathbb{R}^n \times \mathbb{R}^r \times \mathbb{R} \to \mathbb{R}^n$ is a vector field, $\Omega \subset \mathbb{R}^r$ is a bounded set of admissible controls and $x_0 \in \mathbb{R}^n$ is a vector of initial conditions. Suppose that behavior of the system is considered only in the time range $t \in [t_0, t_f]$ where $t_0$ is the initial time and $t_f$ is the final time. If the dynamical system under consideration is non-smooth in some regions of the state space, it may be necessary to supplement the equation (1) with appropriate transition rules, which define evolution of the system in these areas.

Any piecewise smooth function $u: [t_0, t_f] \to \Omega$ is called an admissible control function, or shortly, an admissible control. Assume that for any admissible control $u$ there exists a unique solution $x: \mathbb{R} \to \mathbb{R}^n, t \geq t_0, x(t_0) = x_0$ called a trajectory of the system (1). Obviously, different controls $u$ produce different trajectories $x$. Some of them may be desirable and some may not. In order to assess controls and resulting trajectories qualitatively, it is necessary to propose a performance measure [1]. The performance measure is defined in the following form:

$$J = h[x(t_f), t_f] + \int_{t_0}^{t_f} g[x(t), u(t), t] dt \tag{2}$$

where $h: \mathbb{R}^n \times \mathbb{R} \to \mathbb{R}$ and $g: \mathbb{R}^n \times \mathbb{R}^r \times \mathbb{R} \to \mathbb{R}$ are scalar functions. It is assumed that an admissible control $u$, which produces smaller value of $J$, is better than other admissible control which yields larger $J$. An optimal control $u^*$ is the control which causes that the performance measure $J$ attains a global minimum.

Suppose that a single component $u_j, 1 \leq j \leq r$ of an admissible control $u = [u_1, u_2, \ldots, u_r]^T$ is analyzed. As $u_j$ is piecewise smooth and bounded, it is also integrable [35]. Therefore, it can be approximated in terms of the Fourier series:

$$\tilde{u}_j(t) = \frac{a_{j0}}{2} + \sum_{k=1}^{K} a_{jk} \cos(k\omega t) + \sum_{k=1}^{K} b_{jk} \sin(k\omega t) \tag{3}$$

where $\omega = \frac{2\pi}{t_f - t_0}$, $a_{j0}, a_{jk}, b_{jk} \in \mathbb{R}$, $K \in \mathbb{N}_+$ and $\tilde{u}_j$ is an approximation of $u_j$. If the function $u_j$ satisfies Dirichlet conditions [35], then $\tilde{u}_j$ converges to $u_j$ as $K \to \infty$.[2]

Assume that an optimal control $u^*$ exists. In order to estimate it, the following approach could be utilized.

---

[2] Except, possibly, the points where $u_i^*$ is discontinuous. At such points $t$, the series converges to the value $\lim_{h \to 0^+} \frac{u_i^*(t+h) + u_i^*(t-h)}{2}$.

1) Define each component of the control function $u_j$ in terms of Fourier coefficients $a_{j0}, a_{jk}, b_{jk}, 1 \leq k \leq K$ using the formula (3).
2) Optimize the performance measure $J$ with respect to the parameters $a_{j0}, a_{jk}, b_{jk}, 1 \leq k \leq K$ of each control component $\tilde{u}_j$, i.e. find such values of $a_{j0}, a_{jk}, b_{jk}$, for which $J$ attains the global minimum.

However, such approach brings along a new difficulty: how to select the constraints of optimization, i.e. allowable ranges of the parameters $a_{j0}, a_{jk}, b_{jk}$, so that the resulting approximation (3) meets the control constraints $\widetilde{\boldsymbol{u}}(t) = [\tilde{u}_1(t), \tilde{u}_2(t), \ldots, \tilde{u}_m(t)]^T \in \boldsymbol{\Omega}$ for any $t \in [t_0, t_f]$ ? In order to answer this question, appropriate definitions and theorems should be introduced beforehand.

## 2.2. Mathematical toolbox

Firstly, two features of an admissible control are going to be considered: *shape* and *span*, since they play an important role in the proposed solution of control optimization problem. Moreover, as the presented algorithm involves application of spherical coordinates in an n-dimensional sphere, fundamental information in this topic is going to be provided.

**Definition 1.** *The span of a bounded function $f: \mathbb{D} \to \mathbb{R}$ is the following ordered pair.*

$$(i_f, s_f) = (\inf\{f(\mathbb{D})\}, \sup\{f(\mathbb{D})\}) \qquad (4)$$

The span of a function $f$ is understood as an ordered pair, which contains the infimum and the supremum [35] of its set of values (see Fig. 1). As the image $f(\mathbb{D})$ of the bounded function $f$ is a bounded subset of $\mathbb{R}$, existence and uniqueness of both, $\inf\{f(\mathbb{D})\}$ and $\sup\{f(\mathbb{D})\}$, are guaranteed [35]. Please note that the assumed definition of the span is not connected with the term *linear span*, known from the linear algebra [35].

**Definition 2.** *The shape of a bounded function $f: \mathbb{D} \to \mathbb{R}$, whose span equals $(i_f, s_f)$ and $s_f > i_f$, is the function $\bar{f}: \mathbb{D} \to [0, 1]$ defined as follows.*

$$\bar{f}(x) = \frac{f(x) - i_f}{s_f - i_f} \qquad (5)$$

The assumption $s_f > i_f$ implies that the shape is not defined for constant functions. Intuitively, the shape $\bar{f}$ of the function $f$ is its normalization to the interval $[0, 1]$. The notions of span and shape are illustrated in Fig. 1 below.

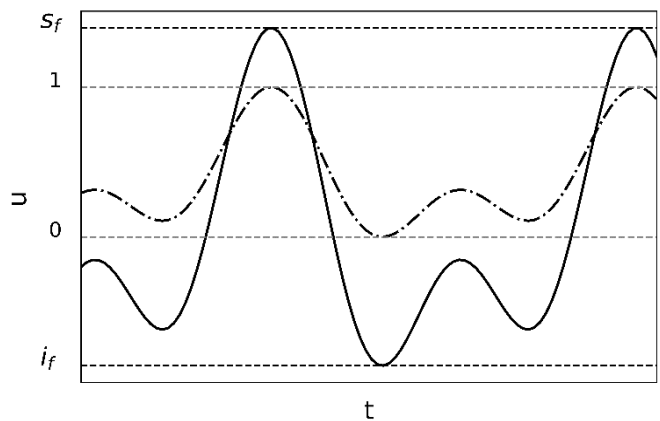

**Fig. 1.** *An exemplary function $f$ (solid line), its shape (dash-dotted line) and span $(i_f, s_f)$*

Note that the set of values of $\bar{f}$, i.e. $\bar{f}(\mathbb{D})$, may not include 0 or 1 if the infimum or the supremum of $f(\mathbb{D})$ is not equal to the minimum or the maximum of $f$ respectively. Moreover, as $f$ can be discontinuous, the set $\bar{f}(\mathbb{D})$ might not be connected [35]. Nevertheless, the property $\bar{f}(\mathbb{D}) \subset [0, 1]$ holds for any bounded, non-constant function $f: \mathbb{D} \to \mathbb{R}$.

Along with the notions of span and shape, it is useful to introduce the following three theorems.

**Theorem 1**. *Any two bounded functions $f: \mathbb{D}_f \to \mathbb{R}, g: \mathbb{D}_g \to \mathbb{R}$ are of the same span if and only if the following relation holds.*

$$\inf\{f(\mathbb{D}_f)\} = \inf\{g(\mathbb{D}_g)\} \wedge \sup\{f(\mathbb{D}_f)\} = \sup\{g(\mathbb{D}_g)\} \tag{6}$$

**Proof.** *The proof of this theorem results directly from Definition 1 and the basic property of the ordered pair [4], which implies that $(i_f, s_f) = (i_g, s_g)$ if and only if $i_f = i_g \wedge s_f = s_g$.*

In other words, two functions are of the same span if both: the infimum and the supremum [35] of their images (sets of values) are equal.

**Theorem 2**. *Any two bounded, non-constant functions $f: \mathbb{D}_f \to \mathbb{R}, g: \mathbb{D}_g \to \mathbb{R}$ are of the same shape if and only if $\mathbb{D}_f = \mathbb{D}_g = \mathbb{D}$ and there exist two constants $a \in \mathbb{R}_+, b \in \mathbb{R}$ such that for any $x \in \mathbb{D}$ the following relation holds.*

$$g(x) = a * f(x) + b \tag{7}$$

**Proof.** *Assume that there exist $a \in \mathbb{R}_+, b \in \mathbb{R}$ such that $g(x) = a * f(x) + b$ for all $x \in \mathbb{D}$. Then, using Definition 2, the following result can be obtained:*

$$\bar{g}(x) = \frac{g(x) - i_g}{s_g - i_g} = \frac{a * f(x) + b - a * i_f - b}{a * s_f + b - a * i_f - b} = \frac{f(x) - i_f}{s_f - i_f} = \bar{f}(x)$$

*as expected. Conversely, assume that $\bar{g}(x) = \bar{f}(x)$. Then, for all $x \in \mathbb{D}$, the following holds:*

$$\frac{g(x) - i_g}{s_g - i_g} = \frac{f(x) - i_f}{s_f - i_f} \to g(x) = \frac{s_g - i_g}{s_f - i_f} f(x) + i_g - i_f \frac{s_g - i_g}{s_f - i_f}$$

*where $\frac{s_g - i_g}{s_f - i_f} \in \mathbb{R}_+$ and $\left(i_g - i_f \frac{s_g - i_g}{s_f - i_f}\right) \in \mathbb{R}$. Thus, the proof is completed.*

The theorem above implies that multiplication of a bounded, non-constant function by a positive constant, as well as adding a constant to it, does not change its shape.

**Theorem 3**. *A bounded, non-constant function $f: \mathbb{D} \to \mathbb{R}$ is uniquely defined by its domain $\mathbb{D}$, its span $(i_f, s_f)$ and its shape $\bar{f}: \mathbb{D} \to [0, 1]$.*

**Proof.** *Suppose that the domain $\mathbb{D}$, the shape $\bar{f}: \mathbb{D} \to [0, 1]$ and the span $(i_f, s_f)$ of a function $f: \mathbb{D} \to \mathbb{R}$ are known. Then, using Definition 2, the following holds for all $x \in \mathbb{D}$:*

$$\bar{f}(x) = \frac{f(x) - i_f}{s_f - i_f} \to f(x) = \bar{f}(x)(s_f - i_f) + i_f \tag{8}$$

*which completes the proof.*

In accordance with the theorem above, in order to define a function, it is enough to specify its domain, shape and span.

Last but not least, in the further parts of paper, the notion of a unit hypersphere (or a unit n-sphere), is going to be necessary. Moreover, spherical coordinates specifying location of a point on the hypersphere will be utilized. Therefore, for completeness of the paper, it seems reasonable to provide the definition of the unit hypersphere and a theorem, which introduces spherical coordinates in the n-sphere. Note that only the most rudimentary information in this topic is provided, for details please refer to the paper [36] and references therein.

**Definition 3**. *A unit n-sphere (hypersphere) $S^n$ is defined as follows [36].*

$$S^n = \{x \in \mathbb{R}^{n+1} : \|x\| = 1\} \tag{9}$$

The unit n-sphere is the set of points in $\mathbb{R}^{n+1}$, whose distance to the origin equals 1. Note that dimension of the hypersphere embedded in the Euclidean space $\mathbb{R}^{n+1}$ equals $n$, as it is enough to specify $n$ parameters in order to select a point in the n-sphere. This fact is illustrated in Fig. 2 and presented in the following theorem.

**Theorem 4**. *For any point $x \in S^n \subset \mathbb{R}^{n+1}$, there exists a unique n-dimensional vector of spherical coordinates $\boldsymbol{\varphi} = [\varphi_1, \varphi_2, \ldots, \varphi_n]^T$, such that $\varphi_1, \varphi_2, \ldots, \varphi_{n-1} \in [0, \pi]$, $\varphi_n \in [0, 2\pi)$ and the following relations hold [36].*

$$x_1 = \cos(\varphi_1) \tag{10a}$$

$$x_2 = \sin(\varphi_1)\cos(\varphi_2) \tag{10b}$$

$$\ldots$$

$$x_n = \sin(\varphi_1) \ldots \sin(\varphi_{n-1})\cos(\varphi_n) \tag{10c}$$

$$x_{n+1} = \sin(\varphi_1) \ldots \sin(\varphi_{n-1})\sin(\varphi_n) \tag{10d}$$

where $x_1, \ldots, x_{n+1}$ are cartesian coordinates of the point $x$. For the proof of Theorem 4, please refer to the paper [36].

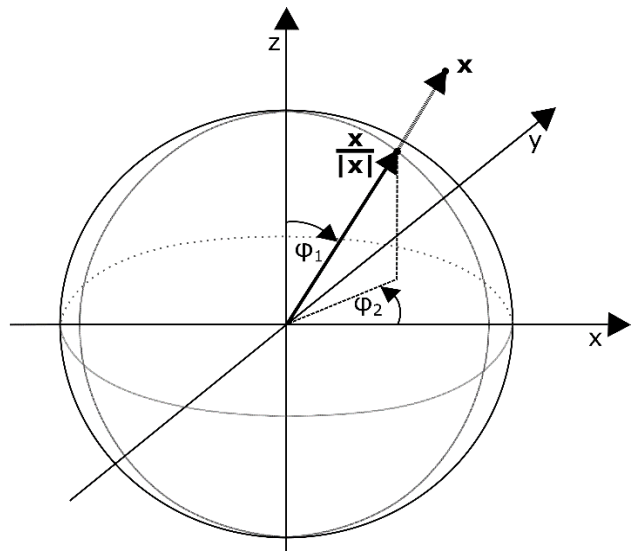

**Fig. 2.** *Illustration of Theorem 4 for the case $n = 2$: the direction of a vector $x \in \mathbb{R}^{n+1} = \mathbb{R}^3$ can be represented by a point in $S^n = S^2$, which is uniquely specified by $n = 2$ spherical coordinates: $\varphi_1 \in [0, \pi], \varphi_2 \in [0, 2\pi)$.*

## 2.3. The algorithm

As the necessary mathematical tools have been introduced, it is time to use them to solve the control optimization problem. According to the Theorem 3, approximation of each component of an admissible control $\tilde{u}_j$ can be uniquely defined by its domain, shape and span. As in this paper control functions are approximated using the Fourier series (3), the domain of such approximation is always $\mathbb{R}$. Therefore, it is enough to investigate the shape and the span. Consequently, if it is possible to parametrize the shape and the span of each component of an admissible control using parameters belonging to well-defined intervals, then searching for the optimal control $\boldsymbol{u}^*$, i.e. minimization of the functional (2), can be reduced to the nonlinear programming problem and can be solved using an appropriate numerical method. Details of such parametrization are provided in this subchapter.

Firstly, the shape of an approximation of the admissible control component $\tilde{u}_j$ is going to be parametrized. Let the formula (3) be transformed to the following form.

$$\tilde{u}_j(t) = \frac{a_{j0}}{2} + [a_{j1}, b_{j1}, \ldots, a_{jK}, b_{jK}][\cos(\omega t), \sin(\omega t), \ldots, \cos(K\omega t), \sin(K\omega t)]^T =$$
$$= \frac{a_{j0}}{2} + \boldsymbol{H}_j [\cos(\omega t), \sin(\omega t), \ldots, \cos(K\omega t), \sin(K\omega t)]^T \quad (11)$$

The attention should be focused at the horizontal vector of amplitudes of subsequent harmonics: $\boldsymbol{H}_j = [a_{j1}, b_{j1}, \ldots, a_{jK}, b_{jK}]$. Assume that the function (11) is non-constant, i.e. $|\boldsymbol{H}_j| > 0$. Let $\overline{\boldsymbol{H}}_j$ be the normalized vector of amplitudes.

$$\overline{\boldsymbol{H}}_j = \frac{\boldsymbol{H}_j}{|\boldsymbol{H}_j|} = [\overline{H}_{j1}, \overline{H}_{j2}, \ldots, \overline{H}_{j(2K-1)}, \overline{H}_{j(2K)}] = \frac{[a_{j1}, b_{j1}, \ldots, a_{jK}, b_{jK}]}{\sqrt{a_{j1}^2 + b_{j1}^2 + \cdots + a_{jK}^2 + b_{jK}^2}} \quad (12)$$

Then, the function:

$$\frac{1}{|\boldsymbol{H}_j|}\left(\tilde{u}_j(t) - \frac{a_{j0}}{2}\right) = \frac{[a_{j1}, b_{j1}, \ldots, a_{jK}, b_{jK}][\cos(\omega t), \sin(\omega t), \ldots, \cos(K\omega t), \sin(K\omega t)]^T}{\sqrt{a_{j1}^2 + b_{j1}^2 + \cdots + a_{jK}^2 + b_{jK}^2}} =$$
$$= \overline{\boldsymbol{H}}_j [\cos(\omega t), \sin(\omega t), \ldots, \cos(K\omega t), \sin(K\omega t)]^T \quad (13)$$

is obtained by adding the constant $\left(-\frac{a_{j0}}{2}\right)$ to $\tilde{u}_j(t)$ (11), followed by multiplying the result by $1/|\boldsymbol{H}_j|$. According to Theorem 2, such operations do not change the shape of the function. Therefore, the function (13) is of the same shape as $\tilde{u}_j(t)$ (11). Consequently, normalization of the vector of harmonics is the operation which preserves the shape of the function. As a result, when geometrical interpretation of the vector $\boldsymbol{H}_j$ is taken into account, it can be stated that the shape of $\tilde{u}_j$ is determined by the *direction* of $\boldsymbol{H}_j$, but not by its *length*. Since the shape of the function (13) depends also on the vector $[\cos(\omega t), \sin(\omega t), \ldots, \cos(K\omega t), \sin(K\omega t)]$, it is influenced by the parameters $\omega$ and $K$. Apart from the direction of $\boldsymbol{H}_j$ and parameters $\omega, K$, there are no other quantities which affect the shape of the function under consideration. Therefore, the shape of (13), equal to the shape of $\tilde{u}_j$ (11), is uniquely defined by the direction of $\boldsymbol{H}_j$ (or, equivalently, by $\overline{\boldsymbol{H}}_j$) and the parameters $\omega, K$.

The normalized vector of harmonics amplitudes $\overline{\boldsymbol{H}}_j$ belongs to the Euclidean space $\mathbb{R}^{2K}$ and indicates a point in this space, whose distance to the origin equals 1. Consider a unit hypersphere $S^{2K-1}$ embedded in $\mathbb{R}^{2K}$, i.e. the set containing these points in $\mathbb{R}^{2K}$, whose distance to the origin equals 1 (see Definition 3). Obviously, the point indicated by $\overline{\boldsymbol{H}}_j$ belongs to this hypersphere, i.e. $(\overline{H}_{j1}, \overline{H}_{j2}, \ldots, \overline{H}_{j(2K-1)}, \overline{H}_{j(2K)}) \in S^{2K-1}$. According to Theorem 4, location of any point on the unit hypersphere $S^{2K-1}$ can be uniquely defined using $2K-1$ spherical coordinates $\boldsymbol{\varphi_j} = [\varphi_{j1}, \varphi_{j2}, \ldots, \varphi_{j(2K-1)}]^T$, such that $\varphi_{j1}, \varphi_{j2}, \ldots, \varphi_{j(2K-2)} \in [0, \pi]$, $\varphi_{j(2K-1)} \in [0, 2\pi)$ and the following relations hold.

$$\bar{H}_{j1} = \cos(\varphi_{j1}) \tag{14a}$$

$$\bar{H}_{j2} = \sin(\varphi_{j1})\cos(\varphi_{j2}) \tag{14b}$$

$$\ldots$$

$$\bar{H}_{j(2K-1)} = \sin(\varphi_{j1})\sin(\varphi_{j2})\ldots\sin(\varphi_{j(2K-2)})\cos(\varphi_{j(2K-1)}) \tag{14c}$$

$$\bar{H}_{j(2K)} = \sin(\varphi_{j1})\sin(\varphi_{j2})\ldots\sin(\varphi_{j(2K-2)})\sin(\varphi_{j(2K-1)}) \tag{14d}$$

The equations above show that the direction $\bar{\boldsymbol{H}}_j$ of the vector of amplitudes $\boldsymbol{H}_j$ is defined by the spherical coordinates $\varphi_{j1}, \varphi_{j2}, \ldots, \varphi_{j(2K-2)} \in [0, \pi]$ and $\varphi_{j(2K-1)} \in [0, 2\pi)$. Consequently, the shape of the function $\tilde{u}_j$ (11) is uniquely defined by the angular coordinates $\boldsymbol{\varphi_j} = [\varphi_{j1}, \varphi_{j2}, \ldots, \varphi_{j(2K-1)}]^T$ and the parameters $\omega, K$. Since the value of the parameter $\omega = \frac{2\pi}{t_f - t_0}$ is known from the problem statement and $K$ can be assumed any positive, natural number[3], the spherical coordinates $\boldsymbol{\varphi_j} = [\varphi_{j1}, \varphi_{j2}, \ldots, \varphi_{j(2K-1)}]^T$ are sufficient to specify the shape of $\tilde{u}_j$.

Now, as the method of shape parametrization has been presented, it is time to consider the span of an admissible control. Obviously, any admissible control must satisfy the condition $\boldsymbol{u}(t) \in \boldsymbol{\Omega}$ for all $t \in [t_0, t_f]$. As in real control systems signals cannot attain arbitrarily large values, it is reasonable to suppose that the set of admissible controls $\boldsymbol{\Omega}$ is bounded. In this paper, it is assumed that $\boldsymbol{\Omega}$ is defined using a set of inequalities in the form:

$$m_j \leq u_j(t) \leq M_j \tag{15}$$

where $m_j, M_j \in \mathbb{R}$, $m_j < M_j$ $1 \leq j \leq r$. In other words, each of $r$ components of a control function $\boldsymbol{u}(t)$ is bounded from above and from below by a constant value. Such constraints are applicable in many practical cases, when outputs of a controller are restricted by fixed parameters.

The inequality (15) presents a condition, which must be met by any admissible control. However, if $(i_{u_j^*}, s_{u_j^*})$ is the span of the $j$-th component of the optimal control $u_j^*$, then obviously $m_j \leq i_{u_j^*}$ and $s_{u_j^*} \leq M_j$, but equalities ($m_j = i_{u_j^*}$, $s_{u_j^*} = M_j$) may not be true. Intuitively, it means that the set of values of $u_j^*$ does not necessarily cover the whole allowable range $[m_j, M_j]$. Therefore, in order to effectively search for the optimal function $u_j^*$, it is necessary to parametrize the span of the approximate control function $\tilde{u}_j$ with respect to the boundaries $m_j, M_j$ in such a way that any value $(i_{\tilde{u}_j}, s_{\tilde{u}_j})$, which ensures that the conditions (15) are fulfilled, can be selected. For this purpose, the parameters defined as follows are used.

$$p_j = \frac{s_{\tilde{u}_j} - m_j}{M_j - m_j} \tag{16}$$

---

[3] As indicated in further parts of the paper, the greater $K$, the better accuracy of the method and the longer optimization time.

$$q_j = \frac{s_{\tilde{u}_j} - i_{\tilde{u}_j}}{s_{\tilde{u}_j} - m_j} \quad (17)$$

Note that $p_j, q_j \in (0, 1]$. Both parameters must be greater than 0, otherwise $s_{\tilde{u}_j} = i_{\tilde{u}_j}$ and the approximate control function would have to be constant. If supremum of the set of values of $\tilde{u}_j$ equals the upper boundary $M_j$, then $p_j = 1$. Analogously, if infimum of the set of values of $\tilde{u}_j$ equals the lower boundary $m_j$, then $q_j = 1$. By choosing appropriate values of $p_j, q_j \in (0, 1]$, any span $(i_{\tilde{u}_j}, s_{\tilde{u}_j})$ of the approximate control function, such that $m_j \leq i_{\tilde{u}_j} < s_{\tilde{u}_j} \leq M_j$, can be obtained. The relation between parameters $p_j, q_j$, the constraints $m_j, M_j$ and the span $(i_{\tilde{u}_j}, s_{\tilde{u}_j})$ is illustrated in Fig. 3.

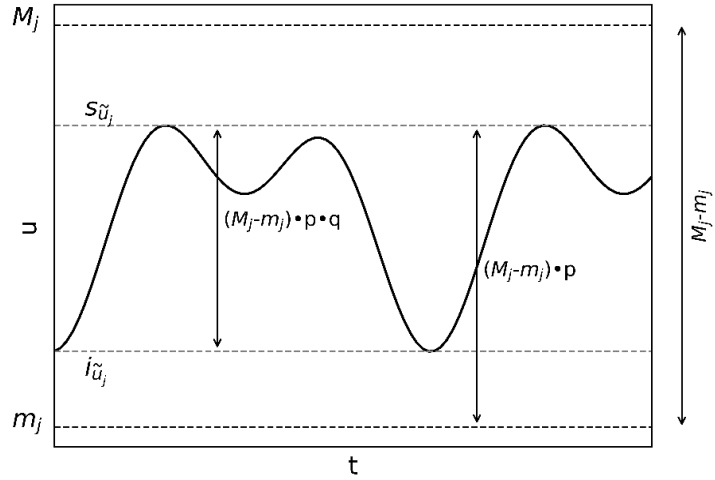

Fig. 3. Illustration of relation between parameters $p_j, q_j$, the constraints $m_j, M_j$ and the span $(i_{\tilde{u}_j}, s_{\tilde{u}_j})$.

The above considerations lead to the following theorem.

**Theorem 5.** Approximation of the j-th component $\tilde{u}_j$ (9) of an admissible control function is uniquely defined by its boundaries $m_j, M_j$, the vector of angular coordinates $\boldsymbol{\varphi_j} = [\varphi_{j1}, \varphi_{j2}, \ldots, \varphi_{j(2K-1)}]^T$ along with the constant $K$ and parameters $p_j, q_j, \omega$.

**Proof.** Suppose that the boundaries $m_j, M_j$, the vector of angular coordinates $\boldsymbol{\varphi_j}$ and the parameters $p_j, q_j, \omega, K$ are known. The shape of $\tilde{u}_j$ (11) equals the shape of (13), which depends exclusively on $\overline{\boldsymbol{H}}_j$ (12) and the parameters $\omega, K$. According to (14a)-(14d), $\overline{\boldsymbol{H}}_j$ is uniquely defined by $\boldsymbol{\varphi_j}$. Consequently, the shape of $\tilde{u}_j$ (11) is determined by $\boldsymbol{\varphi_j} = [\varphi_{j1}, \varphi_{j2}, \ldots, \varphi_{j(2K-1)}]^T$, $\omega$ and $K$. According to (16) and (17), if the boundaries $m_j, M_j$ and parameters $p_j, q_j$ are known, then the span $(i_{\tilde{u}_j}, s_{\tilde{u}_j})$ of $\tilde{u}_j$ can be calculated as follows.

$$s_{\tilde{u}_j} = p_j(M_j - m_j) + m_j \quad (18)$$

$$i_{\tilde{u}_j} = (1 - q_j)s_{\tilde{u}_j} + q_j m_j = (M_j - m_j)(1 - q_j)p_j + m_j \quad (19)$$

Therefore both: the shape and the span of $\tilde{u}_j$ (11) are well defined. Its domain is $\mathbb{R}$. According to Theorem 3, as the shape, the span and the domain of $\tilde{u}_j$ are known, the function $\tilde{u}_j$ is well defined, which completes the proof.

Theorem 5 leads directly to an approximate solution of the control optimization problem. Each component of an admissible control can be approximated by means of the Fourier series (3), (11). The shape of the function under consideration can be found from $\boldsymbol{\varphi}_j, \omega, K$ using formulas (14a)-(14d) and (13). The function (13) is normalized using expression (5), which yields the desired shape[4]. Then, the span is found from parameters $p_j, q_j$ by means of formulas (18), (19). Finally, as the shape and the span are known, $\tilde{u}_j$ is calculated from the expression (8). As a result, definition of the control function approximation $\tilde{u}_j$ in terms of the Fourier parameters $a_{j0}, a_{ji}, b_{ji}, 1 \leq i \leq K$ (3), (11) is obtained. In such a manner, approximation of each component of the control function $\tilde{u}_j$ is determined by $\boldsymbol{\varphi}_j, p_j, q_j, \omega, K$. Consequently, the cost functional $J$ (2) becomes a function of these parameters.

$$J = J(\boldsymbol{\varphi}_j, p_j, q_j, \omega, K), \qquad 1 \leq j \leq r \tag{20}$$

The parameters satisfy the following conditions.

$$\varphi_{jk} \in [0, \pi], \qquad 1 \leq j \leq r, \qquad 1 \leq k \leq 2K - 2 \tag{21a}$$

$$\varphi_{j(2K-1)} \in [0, 2\pi), \qquad 1 \leq j \leq r \tag{21b}$$

$$p_j, q_j \in (0, 1], \qquad 1 \leq j \leq r \tag{21c}$$

$$\omega \geq \frac{2\pi}{t_f - t_0} \tag{21d}$$

$$K \in \mathbb{N}_+ \tag{21e}$$

Note that the value of the fundamental frequency $\omega = \frac{2\pi}{t_f - t_0}$, as it was already explained, enables to approximate any piecewise smooth, bounded function in the time interval $[t_0, t_f]$. However, if the optimal control function turns out to be periodic or almost periodic with period $T \ll t_f - t_0$, then optimization of the parameter $\omega$ with the constraint (21d) can lead to more accurate results under the same number of harmonics $K$, comparing to optimization with fixed $\omega = \frac{2\pi}{t_f - t_0}$. In the case of free parameter $\omega$, its upper limit is not restricted. In practice, the angular frequency of the highest harmonic, i.e. $K\omega$, should not exceed the cutoff frequency of the controller, which is going to physically generate the optimized control function. Otherwise, the controller would not be able to reproduce the desired control accurately.

The parameter $K$ determines the number of harmonics in the control approximation (3), (11). The greater $K$, the better accuracy of the method and the longer optimization time. In practice, it is suggested to start with a small value of $K$ and increase it if necessary.

In order to estimate the optimal control, the cost function (20) must be globally minimized with respect to its parameters $\boldsymbol{\varphi}_j, p_j, q_j, 1 \leq j \leq r$ under constraints (21a)-(21c) with a fixed value $K$. The fundamental frequency $\omega$ can be either fixed ($\omega = \frac{2\pi}{t_f - t_0}$) or optimized under the constraint (21d). Note that the cost function $J$ (20) may not be continuous with respect to its parameters and multiple minima of this function may exist. Therefore, non-gradient, global optimization methods should be used, such as the Differential Evolution [37] or the Particle Swarm Optimization [38].

---

[4] Normalization of the function (11) requires estimation of its span, which can be easily done numerically.

## 3. Numerical example

In this section equations of the control object are derived, details of the simulation are provided and the obtained results are presented.

### 3.1. System equations

In this paper, open loop control of the capsule drive presented in Fig. 4 is optimized. The idea of the system under consideration is based on the papers [28-34]. In this device, swinging of the mathematical pendulum produces inertia forces, which can trigger motion of the whole capsule due to presence of dry friction between the capsule and the underlying surface. Motion of the pendulum can be induced by an external torque $F_\theta(t)$ to be determined.

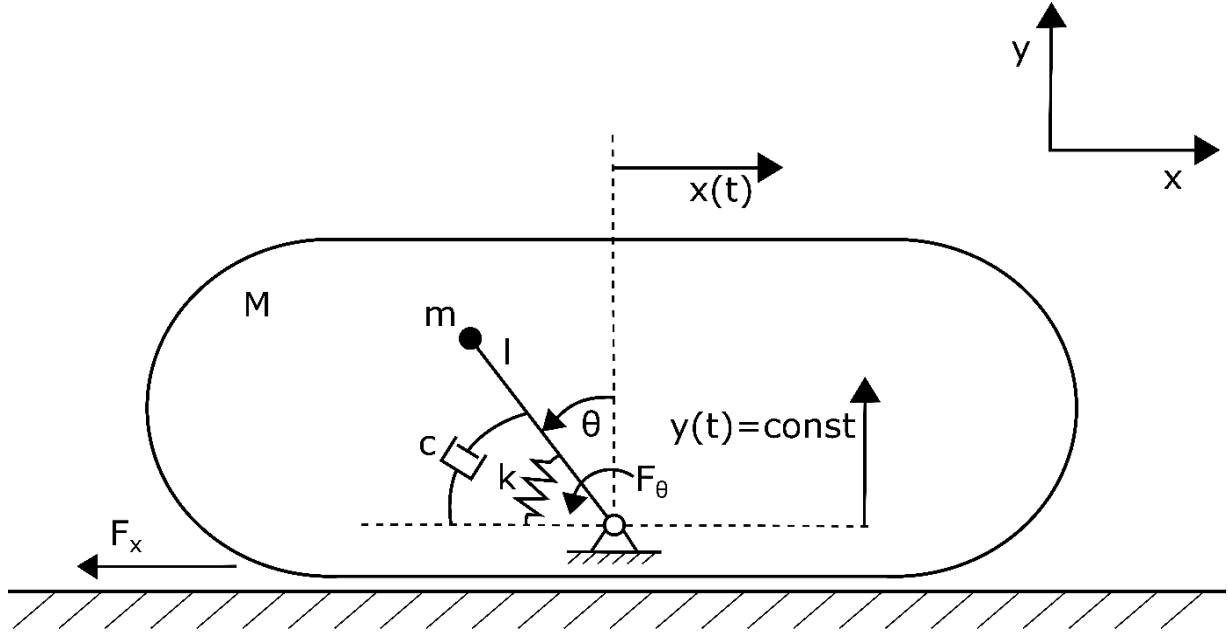

Fig. 4. Scheme of the capsule drive system. M - mass of the capsule, m – mass of the pendulum, l – length of the pendulum, θ – pendulum angle, k – spring stiffness, c – damping coefficient, $F_\theta$ – external torque acting on the pendulum, $F_x$ – friction force, $x(t), y(t)$ – coordinates of the capsule.

Equations of motion of the system are derived using the Lagrange approach and the contact force between the capsule and the ground (the constraint force) is determined by means of the Lagrange multipliers method [6]. Firstly, the Lagrange function is defined as if the vertical coordinate of the capsule $y(t)$ was variable. Then, in the Lagrange equation, the constraint $y(t) = const$ is imposed along with Lagrange multipliers. Two of them vanish from the equations, as there are no constraints connected with the coordinates $x(t)$ and $\theta(t)$, yet the one related to the coordinate $y(t)$ attains the value of the contact load.

The coordinates of the mass center of the pendulum are as follows.

$$x_c(t) = x(t) - l \sin \theta(t), \qquad y_c(t) = y(t) + l \cos \theta(t) \qquad (22)$$

Then, the total kinetic energy of the system can be defined.

$$T = \frac{1}{2} M [\dot{x}^2(t) + \dot{y}^2(t)] + \frac{1}{2} m [\dot{x}_c^{\,2}(t) + \dot{y}_c^{\,2}(t)] \qquad (23)$$

The potential energy of the system is as follows.

$$U = Mgy(t) + mgy_c(t) + \frac{k\theta^2(t)}{2} \tag{24}$$

In the system the following constraint is present.

$$f = y(t) = const \tag{25}$$

By means of the Lagrange function:

$$L = T - U \tag{26}$$

equations of motion and value of the contact force can be obtained by adding non-potential loads $F_x, c\dot{\theta}, F_\theta$ to the appropriate Lagrange equation with constraint forces [6]:

$$\frac{d}{dt}\left(\frac{\partial L}{\partial \dot{\theta}(t)}\right) = \lambda_\theta \frac{\partial f}{\partial \theta(t)} + \frac{\partial L}{\partial \theta(t)} - c\dot{\theta}(t) + F_\theta(t) \tag{27a}$$

$$\frac{d}{dt}\left(\frac{\partial L}{\partial \dot{x}(t)}\right) = \lambda_x \frac{\partial f}{\partial x(t)} + \frac{\partial L}{\partial x(t)} - F_x(t) \tag{27b}$$

$$\frac{d}{dt}\left(\frac{\partial L}{\partial \dot{y}(t)}\right) = \lambda_y \frac{\partial f}{\partial y(t)} + \frac{\partial L}{\partial y(t)} \tag{27c}$$

where $\lambda_\theta, \lambda_x, \lambda_y$ are the Lagrange multipliers. By substitution of expressions (25), (26) to equations (27a)-(27c), taking into account that $y(t) = const$ and $\dot{y}(t) = \ddot{y}(t) = 0$, two equations of motion (28a), (28b) and a value of the constraint force (29) are obtained.

$$ml^2\ddot{\theta}(t) - ml\ddot{x}(t)\cos\theta(t) = mgl\sin\theta(t) - k\theta(t) - c\dot{\theta}(t) + F_\theta(t) \tag{28a}$$

$$(M+m)\ddot{x}(t) - ml\ddot{\theta}(t)\cos\theta(t) + ml\dot{\theta}^2(t)\sin\theta(t) = -F_x(t) \tag{28b}$$

$$R_y(t) = \lambda_y(t) = (M+m)g - ml\ddot{\theta}(t)\sin\theta(t) - ml\dot{\theta}^2(t)\cos\theta(t) \tag{29}$$

Note that $R_y$ is the vertical, reaction (contact) force between the capsule and the underlying surface. In the system of equations under consideration, this is the constraint force and its value is equal to the Lagrange multiplier $\lambda_y$. The complete derivation, conducted using Maxima software, is available in the research data linked to this paper [41].

In the presented system, Coulomb friction model [20] is adopted. From the equation (28a) it can be noticed that the resultant, horizontal force acting on the capsule due to motion of the pendulum is as follows.

$$R_x(t) = ml\ddot{\theta}(t)\cos\theta(t) - ml\dot{\theta}^2(t)\sin\theta(t) \tag{30}$$

Then, the friction model can be formulated.

$$F_x(t) = \begin{cases} \mu R_y(t) sgn[\dot{x}(t)] \leftrightarrow \dot{x}(t) \neq 0 \\ \mu R_y(t) sgn[R_x(t)] \leftrightarrow \dot{x}(t) = 0 \wedge |R_x(t)| \geq \mu R_y(t) \\ R_x(t) \leftrightarrow |R_x(t)| < \mu R_y(t) \end{cases} \tag{31}$$

The following dimensionless variables and parameters are assumed.

$$\Omega = \sqrt{\frac{g}{l}}, \tau = \Omega t, \gamma = \frac{M}{m}, z = \frac{x}{l}, \rho = \frac{k}{m\Omega^2 l^2}, \nu = \frac{c}{m\Omega l^2},$$

$$f_z = \frac{F_x}{m\Omega^2 l}, u_1 = \frac{F_\theta}{m\Omega^2 l^2}, r_z = \frac{R_x}{m\Omega^2 l}, r_y = \frac{R_y}{m\Omega^2 l} \tag{32}$$

All relations between derivatives with respect to the dimensional time $t$ and the dimensionless time $\tau$ are established in the following manner.

$$\dot{x} = \frac{dx}{dt} = \frac{dx}{d\tau}\frac{d\tau}{dt} = \Omega\frac{dx}{d\tau} = \Omega x', \quad \ddot{x} = \frac{d^2x}{dt^2} = \frac{d}{dt}\left(\frac{dx}{dt}\right) = \frac{d}{d\tau}\left(\Omega\frac{dx}{d\tau}\right)\frac{d\tau}{dt} = \Omega^2\frac{d^2x}{d\tau^2} = \Omega^2 x'' \quad (33)$$

Using expressions (32), (33), the set of equations (28a), (28b), can be presented in the dimensionless, matrix form.

$$\begin{bmatrix} 1 & -\cos\theta(\tau) \\ -\cos\theta(\tau) & \gamma+1 \end{bmatrix}\begin{bmatrix} \theta''(\tau) \\ z''(\tau) \end{bmatrix} = \begin{bmatrix} \sin\theta(\tau) - \rho\theta(\tau) - \nu\theta'(\tau) + u_1(\tau) \\ -\theta'^2(\tau)\sin\theta(\tau) - f_z(\tau) \end{bmatrix} \quad (34)$$

In an analogous manner, definitions of the contact load $R_y$ (29), the resultant horizontal force generated by motion of the pendulum $R_x$ (30) and the friction model $F_x$ (31) are transformed.

$$r_y(\tau) = (\gamma+1) - \theta''(\tau)\sin\theta(\tau) - \theta'^2(\tau)\cos\theta(\tau) \quad (35)$$

$$r_z(\tau) = \theta''(\tau)\cos\theta(\tau) - \theta'^2(\tau)\sin\theta(\tau) \quad (36)$$

$$f_z(\tau) = \begin{cases} \mu r_y(\tau)sgn[z'(\tau)] \leftrightarrow z'(\tau) \neq 0 \\ \mu r_y(\tau)sgn[r_z(\tau)] \leftrightarrow z(\tau) = 0 \land |r_z(\tau)| \geq \mu r_y(\tau) \\ r_z(\tau) \leftrightarrow |r_z(\tau)| < \mu r_y(\tau) \end{cases} \quad (37)$$

Details of transformation of the model to the dimensionless form are available in the research data linked to this paper [41]. The inverse of the inertia matrix from Eq. (34), which can be useful in numerical simulations, is also provided in a PDF file. The equations (34)-(37) form a complete model of the capsule system, whose scheme is presented in Fig. 4.

### 3.2. Simulation and optimization details

In this paper, the fundamental question concerning the system (34)-(37) is as follows: what is the forcing function of the pendulum $u_1$, which causes that the capsule covers the maximum distance in a specified interval of dimensionless time ? An approximate answer is provided in the remaining part of this chapter.

First, parameters of the system must be established. In this paper, the following values have been assumed: $\mu = 0.3, \rho = 2.5, \nu = 1.0, \gamma = 10$. In the system, there is only one control function to be optimized: the dimensionless torque acting on the pendulum $u_1$. Obviously, the forcing function $u_1$ cannot attain arbitrary values. Therefore, it has been assumed that, at any time, the condition $u_1(\tau) \in \Omega = [-4, 4]$ must hold. Consequently, boundaries of the control function are: $m_1 = -4, M_1 = 4$. All the simulations have been conducted in the time interval from $\tau_0 = 0$ to $\tau_f = 100$, starting from zero initial conditions: $\theta(0) = \theta'(0) = z(0) = z'(0) = 0$.

The goal of optimization is to find such allowable control $u_1: [\tau_0, \tau_f] \to \Omega = [-4, 4]$, for which the total distance covered by the capsule is maximum. More formally, the performance measure (2) to be minimized can be defined as follows.

$$J = -|z(\tau_f) - z(\tau_0)| \quad (38)$$

The absolute value results from the fact, that the system is symmetric: if a control $u_1$ results in the distance $z(\tau_f)$, then the control $-u_1$ yields $-z(\tau_f)$. Consequently, the crucial issue is to find a distance of a large absolute value, whereas the direction of motion is not so important, as it can be adjusted by changing the sign of the control function.

In accordance with the subchapter 2.3, in particular Theorem 5, approximation of the control function $\tilde{u}_1$ is uniquely defined by the parameters $\boldsymbol{\varphi_1}, p_1, q_1, \omega, K$. The number of harmonics $K$ is fixed in each optimization. Therefore, the performance measure (38) becomes a function $J(\boldsymbol{\varphi_1}, p_1, q_1, \omega)$, which associates each set of parameters $\boldsymbol{\varphi_1}, p_1, q_1, \omega$ (for a fixed $K$) to the negative absolute distance $-|z(\tau_f) - z(\tau_0)|$ covered by the capsule in the specified interval of time. Obviously, values of this function are estimated numerically. First, using spherical coordinates transformation formulas (14a)-(14d), the function (13) is estimated. This function is of a shape defined by the spherical coordinates $\boldsymbol{\varphi_1}$ together with parameters $\omega, K$, but its span requires adjustment. To do so, the current span of the function is estimated numerically: in one period of the function $T = \frac{2\pi}{\omega}$, 1000 points are evenly spaced and the value of the function is computed in each of them using the formula (13). Then, the largest value of the obtained is treated as an approximation of the supremum, and the smallest value is an estimate of the infimum. These results enable normalization of the function to the interval $[0, 1]$ with use of the formula (5), followed by setting the final span defined by the parameters $p_1, q_1$ according to the expression (8). In such a manner, an allowable control $\tilde{u}_1(\tau)$ (11), uniquely specified in terms of the Fourier parameters $a_{10}, a_{1i}, b_{1i}, 1 \leq i \leq K$, is obtained from $\boldsymbol{\varphi_1}, p_1, q_1, \omega$ (for a fixed $K$). Afterwards, this control is used in simulation of the capsule system (34)-(37) in the time interval from $\tau_0 = 0$ to $\tau_f = 100$. As a result, the final distance $z(\tau_f)$ is obtained.

In the conducted numerical research, conversion of the parameters $\boldsymbol{\varphi_1}, p_1, q_1, \omega, K$ to the control function approximation $\tilde{u}_1(\tau)$ (11) is conducted using a *Python 3* script. Then, the Fourier parameters $a_{10}, a_{1i}, b_{1i}, 1 \leq i \leq K$ are conveyed to a dynamically linked library (DLL), which simulates the system (34)-(37). The library has been created in C++ by means of the *CodeBlocks* software with use of the *Boost::Odeint* library for ODE integration and built using the *g++* compiler from the *mingw* package. Simulation of the system is conducted with use of the RK45 [39] variable step integration method, with the absolute tolerance $10^{-9}$ and the relative tolerance $10^{-12}$. Localization of stick/slip transitions is done by means of the bisection method [40]. The Python 3 script for parameters conversion, as well as the full code of the C++ library, are available in the research data linked to this paper [41]. Moreover, a MATLAB simulation of the system (34)-(37) is provided. Finally, optimization of the function $J(\boldsymbol{\varphi_1}, p_1, q_1, \omega)$, equivalent to optimal control approximation, is conducted with use of the Differential Evolution method [37], implemented in the SciPy package for Python 3. The optimization boundaries are specified by expressions (21a)-(21e). Moreover, the upper bound for the parameter $\omega$ equal 10 is used. The optimization code is also included in the research data [41].

### 3.3. Optimization results

Fig. 5 presents the maximal distance covered by the capsule within the dimensionless time interval $0 \leq \tau \leq 100$, depending on the number of harmonics $K$ in the optimization process.

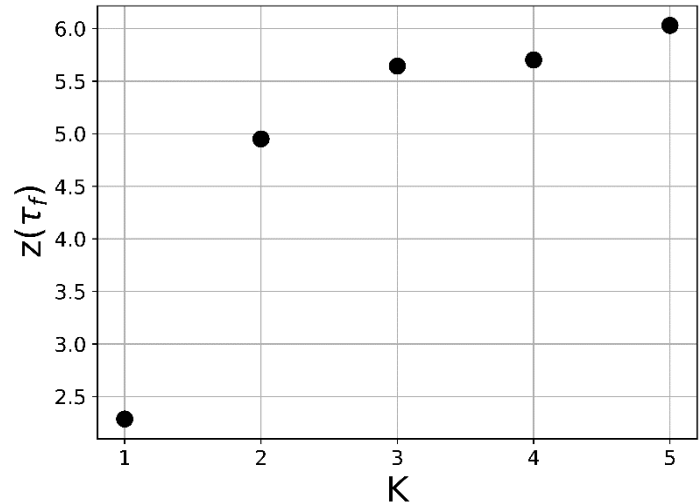

*Fig. 5. Maximal distance covered by the capsule $z(\tau_f)$ vs. number of harmonics used in the optimization process $K$.*

Fig. 5 provides interesting information concerning results of the numerical simulation and, in general, properties of the described algorithm. First of all, the distance covered by the capsule increases with the number of harmonics $K$. This result could be expected from the properties of the method: the larger the number of harmonics in the Fourier approximation $\tilde{u}_1$ (3), (11), the more accurately can the unknown optimal control $u_1^*$ be approximated[5]. Secondly, it can be noticed that the covered distance increases faster for smaller $K$. The result for $K = 2$ is over 100% better then for $K = 1$ (i.e. for the optimal harmonic control), whereas the distance for $K = 5$ is only about 20% longer than for $K = 2$.

Fig. 6 presents how the position of the capsule $z(\tau)$ changes in time, as the control function is optimized using different numbers of harmonics.

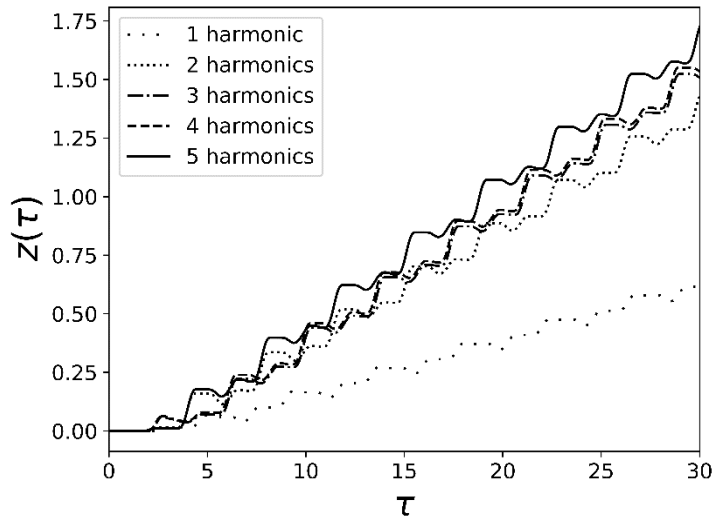

*Fig. 6. Position of the capsule $z(\tau)$ vs. dimensionless time $\tau$.*

---

[5] Obviously, if the optimal control was described exactly in terms of a finite number of harmonics in the Fourier expansion, then the performance measure $J$ would not improve after exceeding this number of harmonics. Otherwise, it is expected that the distance should be increasing (or, equivalently, $J$ should be decreasing) and approaching the optimal value as $K \to \infty$.

Fig. 6 shows that, although the average velocity increases with the number of harmonics $K$, $z(\tau)$ plots have similar features. First of all, none of them is monotonic, which suggests that the optimal control of the capsule drive requires a slight pull back in some moments of time in order to move forward faster later. Moreover, in each plot the stick-slip phenomenon is clearly visible: the flat pieces of graphs (stick) alternate with inclined fragments (slip) in a repetitive manner.

Fig. 7 depicts plots of the optimal control approximations for different values of $K$.

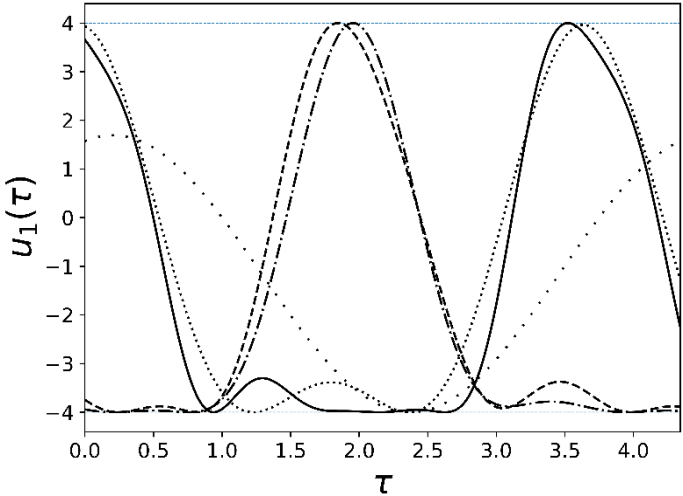

*Fig. 7. Optimal control approximations for different values of K (legend as in Fig. 6).*

The dimensionless time range in Fig. 7 is the longest period from all the estimated controls. Apparently, as $K$ increases, the optimal control approximation deviates from a harmonic function.

## 4. Summary and conclusion

The paper presents a novel, Fourier series based method of open-loop control optimization, which can be used also in non-smooth or discontinuous systems. In the 2$^{nd}$ chapter, the problem is strictly defined, necessary mathematical tools are introduced and the optimization algorithm is explained in details. In the 3$^{rd}$ chapter, a numerical example is presented. The capsule drive, whose motion is triggered by oscillations of an internal mathematical pendulum in presence of the dry friction, is introduced. Its equations of motion are derived using Lagrange method with constraint forces, estimated by means of Lagrange multipliers. Dimensionless parameters are introduced and non-dimensional equations are presented. The simulation procedure is explained and its results are depicted.

The proposed method seems to be a simple and convenient method of open loop optimal control estimation in the presence of control constraints. The presented algorithm has undeniable assets. It is very flexible, very little information about the control object is required. The only requirements are that the set of admissible controls is specified by constant ranges and that for any admissible control it is possible to evaluate the performance measure $J$, whose value is unique. Under these assumptions, it is not even necessary to know equations specifying the control object – it can be a "black box". In particular, non-smooth or even discontinuous systems are not a problem, as long as their solutions exist and are unique for any admissible control. Moreover, the presented method reduces the problem of function optimization to the simpler task of parameters optimization (nonlinear programming), which enables its solution using many global optimization algorithms, such as popular Differential Evolution procedure (implemented in Python's SciPy package). The drawback of the presented algorithm is that properties of the function to be optimized $J(\boldsymbol{\varphi_1}, p_1, q_1, \omega)$, such as its continuity, differentiability, number of minima, etc. are not known. Therefore, it is recommended to use non-gradient global optimization methods, which work relatively slowly and cannot guarantee finding the optimal solution. Nevertheless, the numerical example presented in this paper confirms efficiency of the proposed approach. It is expected that this method will enable approximate investigations of the optimal control in areas, in which it was very difficult or even impossible up to this moment.

## 5. Acknowledgements

This study has been supported by the National Science Centre, Poland, PRELUDIUM Programme (Project No. 2020/37/N/ST8/03448).

This study has been supported by the National Science Centre, Poland under project No. 2017/27/B/ST8/01619.

This paper has been completed while the first author was the Doctoral Candidate in the Interdisciplinary Doctoral School at the Lodz University of Technology, Poland.

# 6. References

[1]: Kirk, D. E. (2004). Optimal control theory: an introduction. Courier Corporation.

[2]: Geering, H. P. (2007). Optimal control with engineering applications. Springer.

[3]: Pontryagin, L. S., Boltyanskii, V. G., Gamkrelidze, R. V., & Mishchenko, E. F. (1962). The mathematical theory of optimal processes, translated by KN Trirogoff. New York.

[4]: Biral, F., Bertolazzi, E., & Bosetti, P. (2016). Notes on numerical methods for solving optimal control problems. IEEJ Journal of Industry Applications, 5(2), 154-166.

[5]: Rao, A. V. (2009). A survey of numerical methods for optimal control. Advances in the Astronautical Sciences, 135(1), 497-528.

[6]: Taylor, J. R. (2005). Classical mechanics. University Science Books.

[7]: Bellman, R. (1954). Dynamic programming and a new formalism in the calculus of variations. Proceedings of the National Academy of Sciences of the United States of America, 40(4), 231.

[8]: Bellman, R., & Kalaba, R. E. (1965). Dynamic programming and modern control theory (Vol. 81). New York: Academic Press.

[9]: Andersson, J. A., Gillis, J., Horn, G., Rawlings, J. B., & Diehl, M. (2019). CasADi: a software framework for nonlinear optimization and optimal control. Mathematical Programming Computation, 11(1), 1-36.

[10]: Leine, R. I., Van Campen, D. H., & Van de Vrande, B. L. (2000). Bifurcations in nonlinear discontinuous systems. Nonlinear dynamics, 23(2), 105-164.

[11]: Stewart, D. E., & Anitescu, M. (2010). Optimal control of systems with discontinuous differential equations. Numerische Mathematik, 114(4), 653-695.

[12]: Clarke, F. H. (1989). Methods of dynamic and nonsmooth optimization (Vol. 57). SIAM.

[13]: Clarke, F. H. (1990). Optimization and nonsmooth analysis (Vol. 5). Siam.

[14]: Frankowska, H. (1984). The maximum principle for a differential inclusion problem. In Analysis and Optimization of Systems (pp. 517-531). Springer, Berlin, Heidelberg.

[15]: Sussmann, H. J. (2004). Optimal control of nonsmooth systems with classically differentiable flow maps. IFAC Proceedings Volumes, 37(13), 543-548.

[16]: Ventura, D., & Martinez, T. (1998, May). Optimal control using a neural/evolutionary hybrid system. In 1998 IEEE International Joint Conference on Neural Networks Proceedings. IEEE World Congress on Computational Intelligence (Cat. No. 98CH36227) (Vol. 2, pp. 1036-1041). IEEE.

[17]: Liu, Y., Wiercigroch, M., Pavlovskaia, E., & Yu, H. (2013). Modelling of a vibro-impact capsule system. International Journal of Mechanical Sciences, 66, 2-11.

[18]: Liu, Y., Pavlovskaia, E., & Wiercigroch, M. (2016). Experimental verification of the vibro-impact capsule model. Nonlinear Dynamics, 83(1), 1029-1041.

[19]: Liu, Y., Jiang, H., Pavlovskaia, E., & Wiercigroch, M. (2017). Experimental investigation of the vibro-impact capsule system. Procedia IUTAM, 22, 237-243.

[20] Liu, Y., Pavlovskaia, E., Hendry, D., & Wiercigroch, M. (2013). Vibro-impact responses of capsule system with various friction models. International Journal of Mechanical Sciences, 72, 39-54.

[21]: Liu, Y., Pavlovskaia, E., Wiercigroch, M., & Peng, Z. (2015). Forward and backward motion control of a vibro-impact capsule system. International Journal of Non-Linear Mechanics, 70, 30-46.

[22]: Nguyen, K. T., La, N. T., Ho, K. T., Ngo, Q. H., Chu, N. H., & Nguyen, V. D. (2021). The effect of friction on the vibro-impact locomotion system: modeling and dynamic response. Meccanica, 1-17.

[23]: Chávez, J. P., Liu, Y., Pavlovskaia, E., & Wiercigroch, M. (2016). Path-following analysis of the dynamical response of a piecewise-linear capsule system. Communications in Nonlinear Science and Numerical Simulation, 37, 102-114.

[24]: Liu, Y., Islam, S., Pavlovskaia, E., & Wiercigroch, M. (2016). Optimization of the vibro-impact capsule system. Strojniški vestnik: journal of mechanical engineering, 62(7-8).

[25]: Maolin, L., Yao, Y., & Yang, L. (2018). Optimization of the Vibro-Impact Capsule System for Promoting Progression Speed. In MATEC Web of Conferences (Vol. 148, p. 10002). EDP Sciences.

[26]: Liu, Y., & Chávez, J. P. (2017). Controlling multistability in a vibro-impact capsule system. Nonlinear Dynamics, 88(2), 1289-1304.

[27]: Liu, Y., Chávez, J. P., Zhang, J., Tian, J., Guo, B., & Prasad, S. (2020). The vibro-impact capsule system in millimetre scale: numerical optimisation and experimental verification. Meccanica, 55(10), 1885-1902.

[28]: Liu, P., Yu, H., & Cang, S. (2018). On the dynamics of a vibro-driven capsule system. Archive of Applied Mechanics, 88(12), 2199-2219.

[29]: Liu, P., Yu, H., & Cang, S. (2019). Modelling and analysis of dynamic frictional interactions of vibro-driven capsule systems with viscoelastic property. European Journal of Mechanics-A/Solids, 74, 16-25.

[30]: Liu, P., Yu, H., & Cang, S. (2016, October). Modelling and dynamic analysis of underactuated capsule systems with friction-induced hysteresis. In 2016 IEEE/RSJ International Conference on Intelligent Robots and Systems (IROS) (pp. 549-554). IEEE.

[31]: Liu, P., Yu, H., & Cang, S. (2018). Optimized adaptive tracking control for an underactuated vibro-driven capsule system. Nonlinear Dynamics, 94(3), 1803-1817.

[32]: Liu, P., Yu, H., & Cang, S. (2018). Trajectory synthesis and optimization of an underactuated microrobotic system with dynamic constraints and couplings. International Journal of Control, Automation and Systems, 16(5), 2373-2383.

[33]: Liu, P., Yu, H., & Cang, S. (2019). Adaptive neural network tracking control for underactuated systems with matched and mismatched disturbances. Nonlinear Dynamics, 98(2), 1447-1464.

[34]: Liu, P., Neumann, G., Fu, Q., Pearson, S., & Yu, H. (2018, October). Energy-efficient design and control of a vibro-driven robot. In 2018 IEEE/RSJ International Conference on Intelligent Robots and Systems (IROS) (pp. 1464-1469). IEEE.

[35]: Bronshtein, I. N., Semendyayev, K. A. F., Musiol, G., & Mühlig, H. (2015). Handbook of mathematics. Springer.

[36]: Blumenson, L. E. (1960). A derivation of n-dimensional spherical coordinates. The American Mathematical Monthly, 67(1), 63-66.

[37]: Storn, R., & Price, K. (1997). Differential evolution–a simple and efficient heuristic for global optimization over continuous spaces. Journal of global optimization, 11(4), 341-359.

[38]: Chen, C. Y., & Ye, F. (2012, May). Particle swarm optimization algorithm and its application to clustering analysis. In 2012 Proceedings of 17th Conference on Electrical Power Distribution (pp. 789-794). IEEE.

[39]: Dormand, J. R., & Prince, P. J. (1980). A family of embedded Runge-Kutta formulae. Journal of computational and applied mathematics, 6(1), 19-26.

[40]: Parker, T. S., & Chua, L. (2012). Practical numerical algorithms for chaotic systems. Springer Science & Business Media.

[41]: Zarychta, S., Sagan, T., Balcerzak, M., Dabrowski, A., Stefanski, A., Kapitaniak, T. (2021). A novel, Fourier series based method of control optimization and its application to a discontinuous capsule drive model - research data, Mendeley Data, V1, doi: 10.17632/4yntbrmtn7.1